\documentclass[reqno]{amsart}
\usepackage{ifthen,latexsym}
\newcommand{\BlackBoxes}{\global\overfullrule5pt}

\BlackBoxes

\newcommand{\Pa}{{\mathbb P}}
\newcommand{\Q}{{\mathbb Q}}
\newcommand{\C}{{\mathbb C}}
\newcommand{\R}{{\mathbb R}}
\newcommand{\N}{{\mathbb N}}

\newcommand{\Fcal}{{\mathcal F}}

\newcommand{\Hcal}{{\mathcal H}}
\newcommand{\Ical}{{\mathcal I}}

\newcommand{\Ucal}{{\mathcal U}}
\newcommand{\Mcal}{{\mathcal M}}

\newcommand{\Ocal}{{\mathcal O}}

\numberwithin{equation}{section}

\newtheorem{proposition}{Proposition}[section]
\newtheorem{lemma}[proposition]{Lemma}
\newtheorem{theorem}[proposition]{Theorem}
\newtheorem{definition}[proposition]{Definition}

\newtheorem{remark}[proposition]{Remark}
\newtheorem{exampleemph}[proposition]{Example}   
\newenvironment{example}{\begin{exampleemph}\begin{upshape}}{\end{upshape}\end{exampleemph}} 

\begin{document}
\title{On Finite-Dimensional Term Structure Models}
\author[Damir Filipovi{\smash{\'c}} \and Josef Teichmann]{Damir Filipovi\'c and Josef Teichmann}
\address{Damir Filipovic, Department of Operations Research and Financial Engineering, Princeton University, Princeton, NJ 08544-5263. Josef Teichmann, Department of Financial and Actuarial Mathematics, Vienna University of Technology, A-1040 Vienna}
\email{filipo@math.ethz.ch, josef.teichmann@fam.tuwien.ac.at}
\thanks{We thank Liuren Wu for his helpful comments on some practical aspects in fixed income markets.}

\date{September 2001 (first draft); \today\;(this draft)}
\begin{abstract}
In this paper we provide the characterization of all
finite-dimen\-sional Heath--Jarrow--Morton models that admit arbitrary
initial yield curves. It is well known that affine term structure
models with time-dependent coefficients (such as the Hull--White
extension of the Vasicek short rate model) perfectly fit any initial
term structure. We find that such affine models are in fact the only
finite-factor term structure models with this property. We also show
that there is usually an invariant singular set of initial yield curves where the affine term structure model becomes time-homogeneous. We also argue
that other than functional dependent volatility structures -- such as local
state dependent volatility structures -- cannot lead to finite-dimensional realizations. Finally, our geometric point of view is
illustrated by several examples.
\end{abstract}
\maketitle

\section{Introduction}
In this paper we provide the characterization of all finite-dimensional Heath--Jarrow--Morton (HJM) models that
admit arbitrary initial yield curves. This is an extension and completion of a series of results obtained by
Bj\"ork et al.~\cite{bjo/chr:99,bjo/sve:01,bjo/lan:00}, and \cite{fil:invariant,fil:lnm01,fil/tei:01,zab:00}. It is well known that affine term
structure models with time-dependent coefficients (such as the Hull--White extension of the Vasicek short rate model \cite{hul/whi:90}) perfectly fit any initial term structure. We find that such affine models are in fact the only
finite-factor term structure models with this property, under some weak assumptions on the volatility structure. We
also show that there is usually an invariant singular set of initial yield curves where the affine term structure model
becomes time-homogeneous. This is again well known for the classical Vasicek \cite{vas:77} and Cox--Ingersoll--Ross
(CIR) \cite{cox/ing/ros:85} short rate models, where the set of consistent inital curves is given explicitely by
the model parameters.

Practitioners and academics alike have a vital interest in finite-factor term structure models, and the distinction
of time-homogenous and inhomogeneous ones. According to \cite{hei/wu:01} there are two groups of practitioners
in the fixed income market.

{\emph{Fund managers}} trade on the yield curve (buy and sell swaps at different maturities), trying to make money
out of it. They do not believe that all the interest rate market quotes are ``correct''. Instead, they in general
use a time-homogeneous two- or three-factor model, estimate the model parameters from long time series data, and
then update the state variables (factors) each day to fit the current term structure. Hence the term structure is
considered as a derivative based on more fundamental state variables (factors), such as in an equilibrium model.
The discrepancies between the fitted term structure and the market prices are preceived as potential trading
opportunities. For example, if the fitted curve is above the two year and ten year swap rates, but is below the
five year swap rate. Then one does a butterfly trade: receiving the five year rate (as one thinks it is high) and
delivering the two year and ten year rates (as one thinks they are low compared to the five year rate). After this
trade, one usually needs to wait for six months or longer for the rates to ``reverse'' (as predicted by the model)
so that one can make money. Since this is a long term game, the model parameters must not change every day.
Parameters have to be constant. If a parameter is time-varying, it is a factor and one needs to specify its
dynamics so that one can make corresponding adjustments for the hedging. A state variable (factor) is time-varying,
but since one has a stochastic model for its evolution, one can check on a daily basis whether its realized value
lies within a statistical confidence interval or not.

{\emph{Interest rate option traders}}, on the other hand, often take the quoted yield curve data, with
minimal or no smoothing, as model input. To fit the observed yield curve perfectly, they allow some of the
model parameters to be time-inhomogeneous. They intend to hedge away instantly all the risks on the yield curve and
only worry about the risk in the implied volatility structure. Yet, low-dimensionality of the model is desirable, since
the number of factors usually equals the number of instruments one needs to hedge in the model. And the daily adjustment
of a large number of instruments becomes infeasible in practice due to transaction costs. Of course, the model factors have
to represent tradable values. But this can usually be achieved by a coordinate transformation.

An HJM model for the forward curve, $x\mapsto r_t(x)$, is determined by the volatility structure, $r_t \mapsto
(x\mapsto \sigma(r_t,x))$, and the market price of risk. Here $r_t(x)$ denotes the forward rate at time $t$ for
date $t+x$ (this is the Musiela \cite{mus:93} parameterization). That is, the price at time $t$ of a zero-coupon
bond maturing at date $T\ge t$ is given by
\[ P(t,T)=e^{-\int_0^{T-t}r_t(x)\,dx}.\]
It is shown in \cite{fil:lnm01} that essentially every HJM model can be realized as a stochastic equation
\begin{equation}\label{hjm1}
\left\{
\begin{aligned}
dr_t&=\left(\frac{d}{dx} r_t+\alpha_{HJM}(r_t)\right)\,dt+\sigma(r_t)\,dW_t\\
r_0&=r^\ast,
\end{aligned}
\right.
\end{equation}
in a Hilbert space $H_w$ of forward curves. We shall recapture the precise setup below in Section~\ref{sechjm}. The
solution process $(r_t)$ in general cannot be realized by a finite-dimensional state process. An HJM model is said
to admit a {\emph{finite-dimensional realization (FDR)}} at the initial forward curve $r^\ast$ if, roughly
speaking, there exists an $m$-dimensional diffusion state process $Z$ (factors), for some $m\in\N$, and a map
$\phi:\R^m\to H$ such that $r_t=\phi(Z_t)$. Notice that $m$, $Z$ and $\phi$ depend on $r^\ast$. In fact,
we shall be interested in those HJM models that admit an FDR of the same dimension at every initial curve, and this
dimension is minimal in some sense.

One of the most basic examples for a finite-dimensional HJM model that fits any initial curve is the
Hull--White extended Vasicek model \cite{hul/whi:90} for the short rate $R_t:=r_t(0)$,
\[ dR_t=\left(b(t)-\beta R_t\right)\,dt+\rho\,dW_t, \quad R_0=r^\ast(0).\]
Here $W$ is a one-dimensional Wiener process defined on
some stochastic basis $(\Omega,(\Fcal_t),\Fcal,\Q)$,
where $\Q$ is the risk-neutral measure. This example will be completely recaptured in Section~\ref{subsecHWV}, where also all the subsequent functions are explicitly given. The coefficients $\beta\ge 0$ and $\rho>0$ are constant,
and $b(t)$ is a function determined by the initial forward curve~$r^\ast$.
The corresponding volatility structure is constant,
\begin{equation}\label{vasvol}
\sigma(r,x)\equiv\sigma(x)=\rho e^{-\beta x},
\end{equation}
and the forward curve has an affine dependence on $R_t$,
\begin{equation}\label{vas1}
r_t(x)=A_{HWV}(t,x)+B_{HWV}(x)R_t.
\end{equation}
Whence a one-factor affine term structure model with time-dependent coefficients (strictly speaking, the
realization is given by the two-dimensional process $(t,R_t)$). For those inital curves for which $b(t)\equiv b$ we obtain a time-homogeneous process $(R_t)$ and $A_{HWV}(t,x)\equiv A_{HWV}(x)$. In this example we
also can easily perform a coordinate transformation that leads to a different state process $Z$. Indeed (as in
\cite[Section~4.1]{bjo/lan:00}) let
\[ dZ_t=-\beta Z_t\,dt+\rho\,dW_t,\quad Z_0=0.\]
Then we have
\[ R_t=Z_t+e^{-\beta t}r^\ast(0)+\int_0^t e^{-\beta(t-s)}b(s)\,dt,\]
and hence
\begin{equation}\label{vas2}
r_t(x)=\tilde{A}_{HWV}(t,x)+B_{HWV}(x)Z_t,
\end{equation}
for some function $\tilde{A}_{HWV}$. Thus we obtained a simpler state-process $Z$ with a similar (affine) functional form of $r_t(x)$ as in \eqref{vas1}.
But $Z$ is not a tradable value and hence cannot be directly used for hedging.

There is a substantial literature providing sufficient conditions for the existence of finite-dimensional HJM
models (see e.g.~\cite{bjo/sve:01,bjo/lan:00} for further reference). A systematic study from a geometric point of view has been made by Bj\"ork et al.~\cite{bjo/chr:99,bjo/sve:01,bjo/lan:00}, and
\cite{fil:invariant,fil:lnm01,fil/tei:01}, see also \cite{zab:00}. In \cite{bjo/sve:01} Bj\"ork and Svensson give
necessary and sufficient conditions for the existence of FDRs. Their key argument is the classical Frobenius
theorem. Therefore they define a Hilbert space, $\Hcal$, on which $d/dx$ is a bounded linear operator. This space
consists solely of entire analytic functions. It is well known however that the forward curves implied by a CIR
short rate model are of the form $r_t(x)=A_{CIR}(x)+B_{CIR}(x)R_t$, where
\begin{equation}\label{cir1}
A_{CIR}(x)=d\frac{e^{ax}-1}{e^{ax}+c}\quad\text{and}\quad B_{CIR}(x)=\frac{be^{ax}}{(e^{ax}+c)^2},
\end{equation}
for some $a,b,c>0$ and $d\ge 0$ (as will be recaptured in Section~\ref{subsecHWCIR}). But $A_{CIR}$ and $B_{CIR}$ are not
entire analytic functions since their extensions on $\C$ have a singularity at $x=(\ln c+i\pi)/a$. Hence the
CIR forward curves do not belong to $\Hcal$, whence the Bj\"ork--Svensson \cite{bjo/sve:01} setting is too narrow
for the HJM framework.

In \cite{fil/tei:01} we recently succeeded to overcome this difficulty. The aim of the present article is to bring
these results to the attention of the financial community. We shall make clear why, under fairly general
assumptions, finite-dimensional HJM models are necessarily affine, and provide a deeper understanding of some
important examples.

The remainder of the paper is as follows. In Sections~\ref{sechjm} and \ref{secfdr} we give the precise setup and definition for a (local) HJM model and an FDR around an initial curve, respectively. In Section~\ref{fun-sec} we recapture the main results (without proofs) from \cite{fil/tei:01} and \cite{fil/tei:01b}. We only consider functional dependent volatility structures, which includes essentially all interesting examples. This is also justified by Section~\ref{secloc}, where we show that a local state dependent volatility structure cannot admit an FDR, unless it is constant. In Section~\ref{secapp} we illustrate our results by explicit calculations for two-dimensional HJM models, which turn out to be either of Hull--White extended Vasicek or CIR type.

\section{HJM Models}\label{sechjm}
Let $(\Omega,(\Fcal_t),\Fcal,\Q)$ be a filtered probability space satisfying the usual conditions, and $W=(W^1,\dots,W^d)$ a $d$-dimensional Wiener process, $d\ge 1$. We suppose that $\Q$ is the risk neutral measure for the subsequent models. Indeed, for the existence of an FDR it is irrelevant whether we are under the physical measure $\Pa$ or under $\Q\sim\Pa$ (see also \cite[Remark~7.1.1]{fil:lnm01}).

We follow \cite{fil:lnm01}, where $r_t$ is regarded as an element in the Hilbert space of forward curves $H_w$. This space consists of absolutely continuous functions $h:\R_{\ge 0}\to \R$, and is equipped with the norm
\[ \|h\|^2_w:=|h(0)|^2+\int_{\R_{\ge 0}}|\partial_x h(x)|^2w(x)\,dx ,\]
where $w:\R_{\ge 0}\to [1,\infty)$ is a non-decreasing $C^1$-function such that $w^{-1/3}$ is integrable on $\R_{\ge 0}$. For example, $w(x)=e^{\alpha x}$ or $w(x)=(1+x)^\alpha$, for $\alpha>0$ or $\alpha>3$, respectively. The closed subspace $H_{w,0}$ of $H_w$, consisting of functions $h\in H_w$ with $\lim_{x\to\infty}h(x)=0$, has the property that
\[ (g,h)\mapsto g\int h,\]
where $\int h$ denotes the definite integral $x\mapsto \int_0^x h(y)\,dy$, is a continuous bilinear operator from $H_{w,0}\times H_{w,0}$ into $H_w$. Hence whenever $\sigma=(\sigma_1,\dots,\sigma_d)$ is a (bounded and) locally Lipschitz continuous map from $H_w$ into $H_{w,0}^d$ then the HJM drift term,
\begin{equation}\label{alphaHJM}
\alpha_{HJM}(h)=\sum_{i=1}^d \sigma_i(h)\int \sigma_i(h),
\end{equation}
is a well-defined (bounded and) locally Lipschitz continuous map from $H_w$ into $H_w$.

On $H_w$ the (unbounded) closed operator $A=d/dx$ generates the strongly continuous semigroup $(S_t)$ of right-shifts
\[ S_t h=h(\cdot+t),\quad t\ge 0.\]
The same holds for the restriction $A_0$ of $A$ to $H_{w,0}$. It can be shown that the domain of $A$ is $D(A)=\{h\in H_w\cap C^1(\R_{\ge 0})\mid \partial_x h\in H_w\}$ and similarly $D(A_0)=\{h\in H_{w,0}\cap C^1(\R_{\ge 0})\mid \partial_x h\in H_{w,0}\}$. We consider $A^2=A\circ A$ and all higher order powers of $A$ and $A_0$. By induction, we have
\begin{align*}
D(A^\infty)&:=\bigcap_{n\ge 0}D(A^n)=\{h\in H_w\cap C^\infty(\R_{\ge 0})\mid (\partial_x)^n h\in H_w,\,\forall n\ge 1\},\\
D(A_0^\infty)&:=\bigcap_{n\ge 0}D(A_0^n)=\{h\in H_{w,0}\cap C^\infty(\R_{\ge 0})\mid (\partial_x)^n h\in H_{w,0},\,\forall n\ge 1\}.
\end{align*}
It is clear that $D(A^\infty_0)\subset D(A^\infty)$. Equipped with the sequence of seminorms $p_n(h)=\sum_{i=0}^n\|A^i h\|_w$, $n\ge 0$, $D(A^\infty)$ and $D(A^\infty_0)$ become Fr\'echet spaces, and $A$ acts as a bounded linear operator on them. We refer to \cite{eng/nag:00} for the theoretical background. In fact, $D(A^\infty)$ and $D(A^\infty_0)$ are in some sense the largest subspaces of $H_w$ and $H_{w,0}$, respectively, with this property. The functions from \eqref{vas1}, \eqref{vas2} and \eqref{cir1}, $A_{HWV}(t,\cdot)$, $\tilde{A}_{HWV}(t,\cdot)$, $A_{CIR}$, $B_{HWV}$, $B_{CIR}$, all lie in $D(A^\infty)$, the last two even in $D(A^\infty_0)$.

By a {\emph{solution}} $(r_t)$ of \eqref{hjm1} we shall always mean a {\emph{continuous mild solution}} (see \cite{dap/zab:92,fil:lnm01}). That is, an $H_w$-valued continuous adapted process $(r_t)$ which satisfies
\begin{equation}\label{mildsol}
r_t=S_t r^\ast+\int_0^t S_{t-s}\alpha_{HJM}(r_s)\,ds+\sum_{i=1}^d\int_0^tS_{t-s}\sigma_i(r_s)\,dW^i_s.
\end{equation}
In the classical HJM notation, where $f(t,T)=r_t(T-t)$, this becomes the familiar expression
\[ f(t,T)=f(0,T)+\int_0^t \tilde{\alpha}_{HJM}(s,T)\,ds+\sum_{i=1}^d\int_0^t\tilde{\sigma}_i(s,T)\,dW^i_s,\quad 0\le t\le T,\]
where $\tilde{\alpha}_{HJM}(s,T):=\alpha_{HJM}(r_s,T-s)$ and $\tilde{\sigma}_i(s,T):=\sigma_i(r_s,T-s)$. If \eqref{mildsol} only holds locally for $t$ replaced by $t\wedge\tau$, for some stopping time $\tau>0$, then $(r_t)$ is a {\emph{local solution}}.

Since the volatility structure $\sigma=(\sigma_1,\dots,\sigma_d)$ determines $\alpha_{HJM}$ by \eqref{alphaHJM}, the following terminology is justified.

\begin{definition}\label{defhjm}
Let $\Ucal$ be a convex open set in $H_w$. A {\emph{(local) HJM model in $\Ucal$}} is a map $\sigma=(\sigma_1,\dots,\sigma_d):\Ucal\to H_{w,0}^d$ such that \eqref{hjm1} admits a unique $\Ucal$-valued (local) solution for every initial curve $r^\ast\in \Ucal$.
\end{definition}

It is shown in \cite[Section~5.2]{fil:lnm01} that $\sigma:\Ucal\to H_{w,0}^d$ is a local HJM model in $\Ucal$ if it is locally Lipschitz continuous. The restriction to $\Ucal$ is convenient since it allows to incorporate important examples such as the CIR model, where $\Ucal$ is the half space $\{h\in H_w\mid h(0)>0\}$.

\section{Finite-Dimensional Realizations}\label{secfdr}
In this section we give the rigorous definition of an FDR as it was sketched in the introduction.

Let $\Ucal$ be a convex open set in $H_w$ and $\sigma$ a local HJM model in $\Ucal$. Let $r^\ast_0\in \Ucal\cap D(A^\infty)$ and $n\in\N$.
\begin{definition}\label{ndr}
We say that $\sigma$ admits an {\emph{$n$-dimensional realization around $r^\ast_0$}} if there exists an open neighborhood $V$ of $r^\ast_0$ in $\Ucal\cap D(A^\infty)$, an open set $U$ in $\R_{\ge 0}\times\R^{n-1} $, and a $C^\infty$-map $\alpha:U\times V\to \Ucal\cap D(A^\infty)$ such that
\begin{enumerate}
\item $r\in \alpha(U,r)$ for all $r\in V$,
\item $D_z\alpha(z,r):\R^n\to D(A^\infty)$ is injective for every $(z,r)\in U\times V$,
\item $\alpha(z_1,r_1)=\alpha(z_2,r_2)$ implies $D_z\alpha(z_1,r_1)(\R^n)=D_z\alpha(z_2,r_2)(\R^n)$ for all $(z_i,r_i)\in U\times V$,
\item for every $r^\ast\in V$ there exists a $U$-valued diffusion process $Z$ and a stopping time $\tau>0$ such that
\begin{equation}\label{rinM}
r_{t\wedge\tau}=\alpha(Z_{t\wedge\tau},r^\ast)
\end{equation}
is the (unique) local solution of \eqref{hjm1} with $r_0=r^\ast$.
\end{enumerate}
\end{definition}
Thus we only consider FDRs in $\Ucal\cap D(A^\infty)$. This seems to be a restriction first, since the original local HJM model is defined on $\Ucal$. However, we shall see in Proposition~\ref{prop2} below that in most interesting cases the FDRs are necessarily in $\Ucal\cap D(A^\infty)$. Also we should not worry about $U\subset\R_{\ge 0}\times\R^{n-1}$. It will become clear that the first component of the diffusion $Z$ can always be chosen to be the time $t\ge 0$.

Definition~\ref{ndr} states that an $n$-dimensional realization around $r^\ast_0$ implies the existence of an FDR at every point $r^\ast$ in a neighborhood of $r^\ast_0$, and these FDRs have a smooth dependence on $r^\ast$. In fact, by i) and ii), each $\alpha(\cdot,r^\ast):U\to \Ucal\cap D(A^\infty)$ is (after a localization) the parametrization of an $n$-dimensional submanifold with boundary, say $\Mcal_{r^\ast}$, of $\Ucal\cap D(A^\infty)$, and \eqref{rinM} says that
\[ r_{t\wedge\tau}\in \Mcal_{r^\ast},\quad\forall t\ge 0.\]
Condition iii) implies that two such leafs $\Mcal_{r_1}$ and $\Mcal_{r_2}$ can only intersect at points where their tangent spaces coincide. According to \cite{fil/tei:01}, the family $\{\Mcal_{r}\}_{r\in V}$ is called an $n$-dimensional {\emph{weak foliation}} on~$V$.

The existence of FDRs around a point is assured by an extended version of the Frobenius theorem (\cite{fil/tei:01}) on the Fr\'echet space $D(A^\infty)$. The Frobenius theorem has also been used by Bj\"ork et al.~\cite{bjo/sve:01,bjo/lan:00} on the Hilbert space $\Hcal$, which however has the drawbacks mentioned in the introduction.

A striking feature of the Frobenius theorem is that it brings together an algebraic condition (dimension of a Lie-algebra) with the analytic problem of the existence of an FDR (weak foliation). The former condition can in many cases be explicitly checked, as it is exemplarily carried out in Section~\ref{secloc}. We do not intend to go further into the theory of Lie-algebras here, but refer the interested reader to \cite{fil/tei:01}. We only mention that, under the appropriate assumptions (in particular, $\sigma:\Ucal\cap D(A^\infty)\to D(A_0^\infty)^d$ has to be $C^\infty$), a {\emph{necessary}} condition for the existence of an $n$-dimensional realization around $r^\ast_0$ is that the Lie-algebra generated by the vector fields
\begin{equation}\label{mudef}
\mu(r):=A r+\alpha_{HJM}(r)-\frac{1}{2}\sum_{i=1}^d D\sigma_i(r)\sigma_i(r)
\end{equation}
and $\sigma_1,\dots,\sigma_d$ on $\Ucal\cap D(A^\infty)$ is at most $n$-dimensional. This algebraic condition typically yields an obstruction for the examined model to admit FDRs, as it is shown in Section~\ref{secloc}.

\section{Functional Dependent Volatility}\label{fun-sec}
In this section we have the idea of $\sigma$ being sensitive with respect to linear functionals of the forward curve. That is, we let $\sigma:\Ucal\to H_{w,0}^d$ be a local HJM model in $\Ucal$ and suppose in addition that $\sigma_i(r)=\phi_i(\ell_1(r),\dots,\ell_p(r))$, for some $p\in\N$, where $\phi_i:\Ocal\subset\R^p\to D(A^\infty_0)$ are smooth maps and $\ell_1,\dots,\ell_p$ are continuous linear functionals on $H_w$. This includes essentially all interesting examples. Indeed, we may think of
\[ \ell_i(r)=\frac{1}{x_i}\int_0^{x_i}r(y)\,dy\quad\text{(benchmark yields)}\]
or
\[ \ell_i(r)=r(x_i)\quad\text{(benchmark forward rates)},\]
for some fixed tenor $x_1,\dots,x_p\ge 0$. For a short rate model we simply have $p=1$ and $\ell(r)=r(0)$ (see \eqref{vasvol} for the trivial constant case $\ell(r)\equiv 1$ and $\phi(z)\equiv \rho\exp(-\beta\,\cdot\,)$). This idea is formalized by the following regularity and non-degeneracy assumptions
\begin{description}
\item[(A1)] $\sigma_i=\phi_i\circ\ell$ where $\ell\in L(H_w,\R^p)$, for some $p\in\N$, $\phi_i:\Ocal\to D(A^\infty_0)$ are $C^\infty$-maps, for $1\le i\le d$, and $\Ocal$ is an open set in $\R^p$ containing $\ell(\Ucal)$.
\item[(A2)] $\phi_1(z),\dots,\phi_d(z)$ are linearly independent, for all $z\in\Ocal$.
\item[(A3)] The linear map $(\ell,\ell\circ A,\dots,\ell\circ A^q):D(A^\infty)\to\R^{p(q+1)}$ is open, for every finite $q\ge 0$.
\end{description}
Assumptions {\bf{(A1)}} and {\bf{(A2)}} are clear, only {\bf{(A3)}} needs some further explanation. Intuitively, it says that the following {\emph{interpolation problem}} is well-posed on $D(A^\infty)$: given a smooth curve $g:\R_{\ge 0}\to\R$, for any finite number of data of the form $\Ical=(\ell(g),\ell(\partial_x g),\dots,\ell((\partial_x)^q g))\in\R^{p(q+1)}$ we can find an interpolating function $h\in D(A^\infty)$ with $(\ell(h),\dots,\ell\circ A^q(h))=\Ical$. Hence degenerate examples such as $p=3$ and $\ell(h)=(h(0),h(1),\int_0^1 h(x)\,dx)$ are excluded. Indeed, here we have $\ell\circ A(h)=(\partial_x h(0),\partial_x h(1),h(1)-h(0))$. Thus the rank of $(\ell,\ell\circ A):D(A^\infty)\to\R^6$ is 5, whence this map is not open. Assumption {\bf{(A3)}} is always satisfied for short rate models, since there $p=1$ and $\ell(h)=h(0)$.

\begin{example}
{\emph{Constant direction volatility. }}Let $\lambda_1,\dots,\lambda_d$ be fixed linearly independent vectors in $D(A^\infty_0)$, and
\begin{equation}\label{ddv}
\sigma_i=\phi_i\circ\ell\quad\text{with}\quad\phi_i(z)=\sum_{j=1}^d\Phi_{ij}(z)\lambda_j,\quad z\in \Ocal,
\end{equation}
for some $C^\infty$-map $\Phi=(\Phi_{ij})$ on $\Ocal$ with values in the regular $d\times d$-matrices. This example has been extensively studied by Bj\"ork et al.~\cite{bjo/sve:01,bjo/lan:00}.
\end{example}
The vital importance of this example will be made clear by Proposition~\ref{propddv} and Theorems~\ref{thm1} and \ref{thmloc} below.

Now let $\sigma$ satisfy {\bf{(A1)}}--{\bf{(A3)}}. This implies in particular that $\sigma:\Ucal\to H_{w,0}^d$ is locally Lipschitz continuous and hence a local HJM model in $\Ucal$.

Let $r^\ast_0\in \Ucal\cap D(A^\infty)$ and $n\in\N$. There is a lower bound for the dimension of an FDR around $r^\ast_0$.
\begin{proposition}\label{prop1}
Suppose $\sigma$ admits an $n$-dimensional realization around $r^\ast_0$. Then $n\ge d+1$.
\end{proposition}
\begin{proof}
This is \cite[Proposition~4.8]{fil/tei:01}.
\end{proof}
Hence there exist no one-dimensional realizations around $r^\ast_0$, even if $d=1$. This confirms the well-known fact that a short rate model $(R_t)$ that fits every initial curve $r^\ast$ in a neighborhood of $r^\ast_0$ in $D(A^\infty)$ contains necessarily some time-dependent parameters. Such that $Z_t=(t,R_t)$ yields a 2-dimensional realization around $r^\ast_0$. This also proves a conjecture in \cite{bjo/sve:01} (see Remark~7.1 therein).

Since low-dimensionality of the state process is preferred, it would not make much sense to look for $n$-dimensional realizations with $n>d+1$ (although this can be done under some non-degeneracy assumptions, as it is carried out in \cite{fil/tei:01}). Rather we restrict our attention to {\emph{minimal}} (that is, $(d+1)$-dimensional) realizations around $r^\ast_0$. In this case our focus on FDRs in $D(A^\infty)$ turns out to be no restriction at all.
\begin{proposition}\label{prop2}
Suppose $\Mcal$ is a $(d+1)$-dimensional $C^\infty$-submanifold (with boundary) of $\Ucal$ that is locally invariant for \eqref{hjm1}. That is, for every $r^\ast\in\Mcal$, the $\Ucal$-valued local solution $(r_t)$ of \eqref{hjm1} satisfies
\[ r_{t\wedge\tau}\in\Mcal,\quad\forall t\ge 0,\]
for some stopping time $\tau>0$.
Then necessarily $\Mcal\subset D(A^\infty)$, and $\Mcal$ is a $C^\infty$-submanifold (with boundary) of $D(A^\infty)$.
\end{proposition}
\begin{proof}
This is \cite[Theorem~3.2]{fil/tei:01b}.
\end{proof}

We now can cite \cite[Theorem~4.10]{fil/tei:01}.
\begin{proposition}\label{propddv}
Suppose $\sigma$ admits a $(d+1)$-dimensional realization around $r^\ast_0$. Then $\sigma$ is of the form \eqref{ddv} (constant direction volatility) on $V$, where $V$ is given by Definition~\ref{ndr}.
\end{proposition}

So far our considerations were local. Now let $U$ denote the set of all $r^\ast_0\in \Ucal\cap D(A^\infty)$ around which $\sigma$ admits a $(d+1)$-dimensional realization. Is it possible that $U=\Ucal\cap D(A^\infty)$? In general the answer is no. Indeed, suppose $U$ is connected. Then, by a simple continuity argument, $\sigma$ is of the form \eqref{ddv} with the same $\lambda_1,\dots,\lambda_d$ everywhere on $U$. It then follows from \cite[Theorem~4.10]{fil/tei:01} that $U$ must not intersect with the singular set
\begin{equation}\label{sigma}
\Sigma:=\{h\in \Ucal\cap D(A^\infty)\mid \nu(h)\in\langle\lambda_1,\dots,\lambda_d\rangle\},
\end{equation}
where
\begin{equation}\label{nudef}
\nu(h):=Ah+\alpha_{HJM}(h)
\end{equation}
and $\langle\lambda_1,\dots,\lambda_d\rangle$ denotes the $d$-dimensional subspace spanned by $\lambda_1,\dots,\lambda_d$. The set $\Sigma$ is not empty in general, but it is small.
\begin{lemma}\label{lem1}
The set $\Sigma$ is closed and lies in a subspace $G$ of $D(A^\infty)$ with $(d+1)\le\dim G\le d^2+d+1$.
\end{lemma}
\begin{proof}
This follows from \cite[Lemmas~4.5 and 4.11]{fil/tei:01}.
\end{proof}

We now can summarize and state the main result, which says that only affine term structure models admit a minimal FDR around any initial curve, see \cite[Theorem~4.13]{fil/tei:01}.
\begin{theorem}\label{thm1}
Suppose there exists an open connected set $U$ in $\Ucal\cap D(A^\infty)$ such that $\sigma$ admits a $(d+1)$-dimensional realization around every $r^\ast_0\in U$. Then there exist linearly independent vectors $\lambda_1,\dots,\lambda_d$ in $D(A_0^\infty)$ such that $\sigma$ is of the form \eqref{ddv} on $U$, and $U\cap\Sigma=\emptyset$, where $\Sigma$ is given by \eqref{sigma}. Moreover, for every $r^\ast\in U$ there exists $\epsilon>0$ and a $C^\infty$-map $\Psi(\cdot,r^\ast):[0,\epsilon)\to \Ucal\cap D(A^\infty)$ with $\Psi(0,r^\ast)=r^\ast$, an $\R^d$-valued diffusion process $Z$ with $Z_0=0$ and a stopping time $0<\tau\le\epsilon$ such that
\begin{equation}\label{atsm}
r_{t\wedge\tau}=\Psi(t\wedge\tau,r^\ast)+\sum_{i=1}^d Z^i_{t\wedge\tau}\lambda_i
\end{equation}
is the $\Ucal$-valued local solution of \eqref{hjm1}. Moreover, $\Psi(\cdot,r^\ast)$ is the unique solution to the evolution equation in $D(A^\infty)$, see \eqref{mudef},
\begin{equation}\label{eveq}
\frac{d}{dt} u(t)=\mu(u(t)),\quad u(0)=r^\ast.
\end{equation}

If, in addition, $U=\Ucal\cap D(A^\infty)\setminus\Sigma$ then for every $r^\ast\in\Sigma$ there exists an $\R^d$-valued time-homogeneous diffusion process $Z$ with $Z_0=0$ and a stopping time $0<\tau\le\epsilon$ such that
\begin{equation}\label{atsm2}
r_{t\wedge\tau}=r^\ast(t\wedge\tau)+\sum_{i=1}^d Z^i_{t\wedge\tau}\lambda_i
\end{equation}
is the $\Ucal$-valued local solution of \eqref{hjm1}, and $r_{t\wedge\tau}\in \Sigma$ for all $t\ge 0$. Hence $\Sigma$ is locally invariant for \eqref{hjm1}.
\end{theorem}
Thus here we have the announced decomposition of the space of initial forward curves $\Ucal\cap D(A^\infty)$ into $U$, where the (affine) factor model becomes time-inhomogeneous, and the singular set $\Sigma$, which is invariant for the model dynamics and where the (affine) factor model becomes time-homogeneous. This phenomenon was known for some particular models, and it now is proved in full generality. We shall further illustrate this result for the case $d=1$ in Section~\ref{secapp}.
\begin{remark}\label{rem1}
The vector field $\mu$ in \eqref{eveq} can be replaced by $\nu$ and the statement about the affine form \eqref{atsm} remains true if $\tau$ is chosen such that $r_{t\wedge\tau}$ is $\overline{U}\cap\Ucal$-valued. This follows since $\mu(h)-\nu(h)\in\langle\lambda_1,\dots,\lambda_d\rangle$ for all $h\in\overline{U}\cap\Ucal$ ($=\Ucal$ if $U=\Ucal\setminus\Sigma$), which yields a straightforward modification of $\Psi$ and $Z$ in \eqref{atsm}.
\end{remark}

\section{Local state dependent volatilities}\label{secloc}
Continuing the discussion at the end of Section~\ref{sechjm}, we are formally given smooth vector fields $X_{1},...,X_{m}$ on $ \mathcal{U} \cap
D(A^{\infty})$. Involutivity -- the algebraic condition of the Frobenius Theorem -- amounts to saying that for any
$i,j=1,...,m$ there exist smooth
functions $\lambda_{1},...,\lambda_{m}:\mathcal{U}_0 \to\mathbb{R}$, such that%
\[
\lbrack X_{i},X_{j}](h)=\sum_{l=1}^{m}\lambda_{l}(h)X_{l}(h)
\]
for $h$ in some smaller domain $ \mathcal{U}_0 \subset \mathcal{U} \cap D(A^\infty)$. We recapture some facts, which have been
derived in \cite{fil/tei:01}. The calculation of Lie brackets can be performed on the space
$C^{\infty}(\mathbb{R}_{\geq0},\mathbb{R})$, since $D(A^{\infty}) $ is a subspace of $C^{\infty}(\mathbb{R}_{\geq0},\mathbb{R})$ and the vector fields that we consider are defined on the latter. Lie
brackets are calculated by Fr\'{e}chet
derivatives%
\[
\lbrack X_{i},X_{j}](h)=DX_{i}(h)\cdot X_{j}(h)-DX_{j}(h)\cdot X_{i}(h),
\]
and calculating Fr\'{e}chet derivatives is equivalent to derivation of
$D(A^\infty)$-valued curves%
\[
DX_{i}(h)\cdot v=\frac{d}{dt}|_{t=0}X_{i}(h+tv).
\]

We want to demonstrate that local volatilities cannot lead to finite-dimensional realizations, since the algebraic
conditions are not satisfied. For simplicity we shall assume that $d=1$. Given $\phi:\mathbb{R}%
_{\geq0}\times\mathbb{R}\rightarrow\mathbb{R}$, we define $\sigma (h)(x)=\phi(x,h(x))$ for $h\in \mathcal{U}
\cap D(A^{\infty})$ and $x\geq0$. By cartesian closedness (see \cite{fil/tei:01}) we
obtain for the Fr\'{e}chet derivative%
\begin{align*}
(D\sigma(h)\cdot v)(x)  &  =\frac{d}{dt}|_{t=0}\phi(x,h(x)+tv(x))\\ &
=\phi^{\prime}(x,h(x))v(x)
\end{align*}
for any $x\geq0$. Consequently $D\sigma(h)\cdot v=(\phi^{\prime}\circ
h)v$, where we have an ordinary multiplication of functions. With
$\circ$ we denote the composition with respect to the second variable
and with $.^{\prime}$ we
denote derivation with respect to the second variable. Therefore%
\[
\mu(h)=\frac{d}{dx}h+(\phi\circ h)(\int\phi\circ h)-\frac{1}{2}((\phi )^{\prime}\circ h)(\phi\circ h),
\]
see \eqref{mudef}. Conditions on $ \phi $ can be found such that all parts of the drift vector field and the local volatility are
smooth maps.

\begin{theorem}\label{thmloc}
Assume that $[\mu,\sigma](h)\in\left\langle \sigma(h),\mu(h) \right\rangle $ for $ h $ in some open set $\Ucal_0\subset\mathcal{U} \cap D(A^{\infty})$. Then $ \sigma(h)(x)=\rho \exp(-\beta x) $ for $h\in\Ucal_0$, for some constants $ \rho $ and $ \beta $, which is the Vasicek volatility structure.
\end{theorem}

\begin{proof}
We have to calculate one Lie bracket:%
\begin{align*}
D\mu(h)\cdot v &  =\frac{d}{dx}v+((\phi^{\prime}\circ h)v)(\int_{0}^{.}%
\phi\circ h)+(\phi\circ h)(\int_{0}^{.}(\phi^{\prime}\circ h)v)-\\
&  -\frac{1}{2}(\phi^{\prime\prime}\circ h)(\phi\circ h)\nu-\frac{1}{2}%
(\phi^{\prime}\circ h)^{2}v
\end{align*}
for $h\in \mathcal{U}\cap D(A^\infty)$ and $v\in D(A^{\infty})$. The derivative with respect to the first variable
of $\phi$ is denoted by $\partial_{1}$. This leads to%
\begin{gather*}
\lbrack\mu,\sigma](h)=(\phi^{\prime}\circ h)\frac{d}{dx}h+(\partial_{1}%
\phi\circ h)+(\phi^{\prime}\circ h)(\phi\circ h)\left(  \int\phi\circ
h\right)  +\\ +(\phi\circ h)\left(  \int(\phi^{\prime}\circ
h)(\phi\circ h)\right)
-\frac{1}{2}(\phi^{\prime\prime}\circ h)(\phi\circ h)(\phi\circ h)-\frac{1}%
{2}(\phi^{\prime}\circ h)^{2}(\phi\circ h)-\\ -(\phi^{\prime}\circ
h)\frac{d}{dx}h-(\phi^{\prime}\circ h)(\phi\circ h)\left(
\int(\phi\circ h)\right)  +\frac{1}{2}(\phi^{\prime}\circ
h)(\phi^{\prime}\circ h)(\phi\circ h)\\ =(\partial_{1}\phi\circ
h)+(\phi\circ h)\left(  \int(\phi^{\prime}\circ h)(\phi\circ h)\right)
-\frac{1}{2}(\phi^{\prime\prime}\circ h)(\phi\circ h)^{2}.
\end{gather*}
We now shall evaluate the equation%
\[
\lbrack\mu,\sigma](h)=\lambda(h)\sigma(h)+\lambda_{0}(h)\mu(h)
\]
for $h\in \mathcal{U}_0 $. Assume that $\lambda_{0}(h_{\infty})\neq0$ for some $h_{\infty }\in
\mathcal{U}_0$. We can isolate $\frac{d}{dx}h$, then we take a sequence $(h_{n})_{n\geq0}$ in $\mathcal{U}_0$ with
$h_{n}\rightarrow h_{\infty}$ in $D(A^{m})$ as $n\rightarrow\infty$ but not in
$D(A^{m+1})$, where $m$ is chosen such that $\lambda_{0}(h_{n})\rightarrow\lambda _{0}(h_{\infty})$
(this is possible by applying Theorem 2.3 of \cite{fil/tei:01} and adjusting $\mathcal{U}_0$). Hence the only term
which does not converge up to order $m$ is $\frac{d}{dx}h_{n}$. Consequently $\lambda_{0}=0$ on $\mathcal{U}_0$.
Next we analyze
the resulting equation%
\[
\lbrack\mu,\sigma](h)=\lambda(h)\sigma(h),
\]
where we proceed by the same reasoning.

First we assume that we are given a point $h_0\in \mathcal{U}_0$ and $x_{0}\geq0$ such that
$\phi(x_{0},h_0(x_{0}))\neq0$, then we can divide by $\sigma(h)(x)$ for $h$ in a neighborhood of $h_0$ and $x$ in a neighborhood of $x_{0}$ and identify a
logarithmic derivative. To this
equation we apply again the operator $\frac{d}{dx}$ and obtain%
\begin{gather*}
(\partial_{1}^{2}\ln|\phi|\circ
h)+(\partial_{1}(\ln|\phi|)^{\prime}\circ
h)h^{\prime}+(\phi^{\prime}\circ h)(\phi\circ h)-\frac{1}{2}(\partial_{1}%
\phi^{\prime\prime}\circ h)(\phi\circ h)-\\
-\frac{1}{2}(\phi^{\prime\prime\prime}\circ h)h^{\prime}(\phi\circ
h)-\frac
{1}{2}(\phi^{\prime\prime}\circ h)(\partial_{1}\phi\circ h)-\frac{1}{2}%
(\phi^{\prime\prime}\circ h)(\phi^{\prime}\circ h)h^{\prime}=0.
\end{gather*}
We can again isolate $h^{\prime}=\frac{d}{dx}h$, which leads to a
contradiction as before, therefore its coefficient has to vanish
identically.
This leads to the following two equations:%
\[
(\partial_{1}^{2}\ln|\phi|\circ h)+(\phi\circ h)(\phi^{\prime}\circ
h)-\frac{1}{2}(\partial_{1}\phi^{\prime\prime}\circ h)(\phi\circ
h)-(\phi^{\prime\prime}\circ h)(\partial_{1}\phi\circ h)=0
\]
and%
\[
(\partial_{1}(\ln|\phi|)^{\prime}\circ
h)-\frac{1}{2}(\phi^{\prime\prime \prime}\circ h)(\phi\circ
h)^{2}-(\phi^{\prime\prime}\circ h)(\phi^{\prime }\circ h)=0.
\]
We take these two equations and evaluate them for $h$ in a neighborhood of $ h_0 $ and $x$ in a
neighborhood of $x_{0}$, which leads to the equations%
\begin{align*}
\partial_{1}^{2}\ln|\phi|+\phi\phi^{\prime}-\frac{1}{2}\partial_{1}^{2}%
(\phi^{\prime\prime}\phi)  & =0\\
\partial_{1}(\ln|\phi|)^{\prime}-\frac{1}{2}(\phi^{\prime\prime}\phi
)^{\prime}  & =0.
\end{align*}
Taking the derivative $.^{\prime}$ in the first and $\partial_{1}$ in
the
second and finally the difference of the resulting equations we obtain%
\[
(\phi\phi^{\prime})^{\prime}=0
\]
and therefore $\phi(x,y)\frac{\partial}{\partial y}\phi(x,y)=g(x)$ for $(x,y)$ in a neighborhood of $(x_{0},h_0(x_0))$ with some smooth function $g:\mathbb{R}_{\geq0}\rightarrow\mathbb{R}$. This equation has a smooth
solution if and only if $\frac{\partial}{\partial y}\phi(x,y)=0$. Therefore $\partial_{1}^{2}\ln|\phi|=0$, which
leads to $\phi(x,y)=\rho%
\exp(- \beta x)$ for $(x,y)$ in a neighborhood of $(x_{0},h_0(x_0))$. By continuity we obtain the global result.
\end{proof}

\begin{remark}
The same method can be applied to the $d$-dimensional case. This leads
to a similar assertion, however, the calculus is rather ambitious.
\end{remark}

\section{Characteristic Examples}\label{secapp}

In the seminal papers \cite{bjo/sve:01} and \cite{bjo/lan:00} finite-dimensional realizations, in particular the
Hull-White extensions of the Vasicek and CIR-model, are considered for the first time from the geometric point of view. In
addition to their excellent treatment (compare Section 5 of \cite{bjo/lan:00} or Section 7 of \cite{bjo/sve:01}),
we prove that the Hull-White extensions of the Vasicek and CIR model are the only $2$-dimensional local HJM models and we demonstrate the importance of the corresponding singular sets. The same type of analysis can also be performed in higher dimensional cases, which will be done elsewhere. At the end of this section we provide an example of how to embed the Svensson family as a leaf of a weak foliation associated to a functional dependent volatility structure.

We let again $d=1$. Starting with a functional dependent volatility structure $\sigma$ and $\Ucal$ as in Section \ref{fun-sec} and assuming the existence of
$2$-dimensional realizations on $\Ucal\cap D(A^\infty)\setminus\Sigma$ (see \eqref{sigma}), we necessarily arrive by Proposition \ref{propddv} at a constant direction
volatility on $\Ucal$. We shall show that this volatility is either of the Vasicek or CIR type.

In view of {\bf{(A2)}} we have $\sigma\neq 0$ on $\Ucal$, hence we can write $\sigma(r)=\phi(r)\lambda$ for $r\in\Ucal$, for some $\lambda\in D(A^\infty_0)\setminus\{0\}$ and a smooth map $\phi:\Ucal\to\R$, such that without loss of generality $\phi>0$ (by a slight abuse of notation, the meaning of $\phi$ here is different from Section~\ref{fun-sec}). We want to
specify under which conditions this volatility structure admits $2$-dimensional realizations and how they look
like. This is already done in Section 7.3 of \cite{bjo/sve:01}, however, their special setting does not allow to
treat the CIR-case.

Writing $\psi
(r):=\phi(r)(D\phi(r)\cdot\lambda)$, we obtain  for $ r \in \Ucal\cap D(A^\infty) $ %
\begin{align*}
D\sigma(r)\cdot h  & =(D\phi(r)\cdot h)\lambda\\
D\sigma(r)\cdot\sigma(r)  & =\phi(r)(D\phi(r)\cdot\lambda)\lambda
=\psi(r)\lambda\\
\mu(r)  & =\frac{d}{dx}r+\phi(r)^{2}\lambda\int\lambda-\frac{1}{2}%
\psi(r)\lambda\\ D\mu(r)\cdot h  & =\frac{d}{dx}h+2\phi(r)(D\phi(r)\cdot h)\lambda\int
\lambda-\frac{1}{2}(D\psi(r)\cdot h)\lambda.
\end{align*}
Consequently we can calculate the Lie bracket%
\begin{align*}
\lbrack\mu,\sigma](r)  & =\phi(r)\frac{d}{dx}\lambda+2\phi(r)\psi(r)\lambda
\int\lambda-\frac{1}{2}\phi(r)(D\psi(r)\cdot\lambda)\lambda-\\ &
-(D\phi(r)\cdot\frac{d}{dx}r)\lambda-\phi(r)^{2}(D\phi(r)\cdot\lambda
\int\lambda)\lambda+\frac{1}{2}\psi(r)(D\phi(r)\cdot\lambda)\lambda.
\end{align*}
We assume $[\mu,\sigma](r)\in\langle\lambda\rangle$ on $\mathcal{U}\cap D(A^\infty)$, which follows from the Frobenius condition and is justified by Lemmas 2.12 and 3.4 of
\cite{fil/tei:01}. We can divide by $\phi(r)$ and obtain an equation%
\[
\frac{d}{dx}\lambda+2\psi(r)\lambda\int\lambda-\theta(r)\lambda=0
\]
with some smooth function $ \theta : \mathcal{U}\cap D(A^\infty) \to \mathbb{R} $. There are consequently two cases:
\begin{enumerate}
\item  If $\lambda$ and $\lambda\int\lambda$ are linearly independent in
$D(A^{\infty})$, then by derivation with respect to $r$ we obtain that $\psi$ and $\theta$ are constant, say
$2\psi(r)=a$ and $ \theta(r)=b$ with real numbers $ a $ and $ b $. Defining $\Lambda:=\int\lambda$ we obtain
finally a Riccati equation for $\Lambda$, which yields the CIR-type if
$a\neq0$ or the Vasicek-type if $a=0$:%
\begin{equation}\label{ricc1}
\frac{d}{dx}\Lambda+\frac{a}{2}\Lambda^{2}+b\Lambda=\lambda(0),\quad\Lambda(0)=0.
\end{equation}
The Ho-Lee
model is considered as particular case of the Vasicek model for $ b = 0 $.
\item  If $\lambda$ and $\lambda\int\lambda$ are linearly dependent in
$D(A^{\infty})$, then we necessarily obtain an equation of the type%
\begin{equation*}
\frac{d}{dx}\lambda+b\lambda=0,
\end{equation*}
which yields that $\lambda$ is vanishes identically, since otherwise $\lambda$ and $\lambda\int\lambda$ are
linearly independent. This case was excluded at the beginning.
\end{enumerate}
Notice that by \eqref{ricc1},  $\lambda(0)=0$ if and only if $\lambda=0$, which is not possible. Hence a fortiori we have $\lambda(0)\neq 0$, such that by rescaling we always can assume that $\lambda(0)=1$. This observation slightly improves the discussion in Section 7.3 in \cite{bjo/sve:01}.

By the definition of $ \psi $ we have $D\phi^2(r)\cdot\lambda=a$, hence we obtain the following representation for $ \phi $. We split $ D(A^{\infty}) $ into
$ \R \lambda + E $, where $ E := \ker ev_0 $. We denote by $pr:D(A^\infty)\to E$ the corresponding projection. Then
\begin{equation}\label{phieq}
\phi (h) = \sqrt{aev_0(h) +  \eta(pr(h))} ,
\end{equation}
where $ \eta:
pr(\mathcal{U}\cap D(A^\infty)) \subset E \to \R $ is a smooth function (compare with Proposition 7.3 of \cite{bjo/sve:01}).

Recalling \eqref{sigma}, we have
\[ \Sigma=\left\{h\in\Ucal\cap D(A^\infty)\mid \nu(h)=A h+\phi(h)^2\lambda\int\lambda\in\langle\lambda\rangle\right\}.\]
Thus, if $\lambda$ and $\lambda\int\lambda$ are linearly independent in
$D(A^{\infty})$ then any $h \in \Sigma$ is necessarily of the form%
\[
h=a_{1}+a_{2}\Lambda^{2}+a_{3}\Lambda
\]
in all cases for some real numbers $ a_i $. By the particular representation of $ \phi $ we obtain that $ a_2 =
aa_1 + g(a_3) $, where $ g $ is some smooth real function derived from

$$ aa_1 + \eta(a_2 {\Lambda}^2 + a_3 \Lambda) = a_2. $$

By $ Fl^X $ we denote the local (semi-)flow of a vector field $X$ on $\Ucal\cap D(A^{\infty}) $.
The leaves through $ r^{\ast} $ of the weak foliation are given by the local parametrization%
\[
(u_{0},u_{1})\mapsto Fl_{u_{0}}^{\nu}(r^{\ast})+u_{1}\frac{d}{dx}\Lambda
\]
if $ r^{\ast} $ does not lie in the singular set $ \Sigma $. If $ r^{\ast} \in \Sigma $, then the leaf is a one
dimensional immersed submanifold of $ \langle 1, \Lambda, {\Lambda}^2 \rangle $. Notice that the stochastic
evolution of the factor process takes place in the $u_1$-component, see Theorem~\ref{thm1} and Remark~\ref{rem1}.

We summarize the preceding results in the following theorem.
\begin{theorem}
Let $\sigma$ and $ \mathcal{U}$ be as above. Assume that $\sigma$ admits a 2-dimensional realization around any initial curve $r^\ast\in\Ucal\cap D(A^\infty)\setminus\Sigma$. Then there exists $\lambda\in D(A^\infty_0)$ and a function $\phi:\Ucal\to \R_{>0}$ of the form \eqref{phieq} such that $\sigma(h)=\phi(h)\lambda$. The singular set $\Sigma$ is a (possibly empty) subset of $\langle 1,\Lambda,\Lambda^2\rangle$, where $\Lambda=\int\lambda$ satisfies the Riccati equation \eqref{ricc1}. The local HJM model is an affine short rate model. That is, for every initial curve $r^\ast\in\Ucal\cap D(A^\infty)$ there exist functions $b:\R_{\ge 0}\times\R\to\R$, $\theta:\R_{\ge 0}\to\R$, $A:[0,\epsilon)\times\R_{\ge 0}\to\R$ and a stopping time $\tau>0$ such that
\begin{equation}\label{atsmX}
r_{t\wedge\tau}(x)=A(t\wedge\tau,x)+\lambda(x) R_{t\wedge\tau}
\end{equation}
is the unique $\Ucal$-valued local solution to \eqref{hjm1} and the short rates $R_t=r_t(0)$ follow, locally for $t\le\tau$, a time-inhomogeneous diffusion process
\[ dR_t=b(t,R_t)\,dt+\sqrt{a R_t+\theta(t)}\,dW_t.\]
This process becomes time-homogeneous if and only if $r^\ast\in\Sigma$, and then $r_{t\wedge\tau}\in\Sigma$ for all $t\ge 0$.
\end{theorem}
\begin{proof}
We know that $\lambda(0)\neq 0$. Hence \eqref{atsmX} follows from \eqref{atsm}. The rest of the theorem is a consequence of Theorem~\ref{thm1} and the preceding discussion.
\end{proof}
\subsection{The Hull-White extension of the Vasicek model}\label{subsecHWV}

We consider the vola\-tility structure of the Vasicek model: $\sigma(r)(x)=\rho\exp(-\beta x)=\rho\lambda$ with
$\rho>0$ and $\beta>0$, for $r\in \Ucal\cap D(A^\infty)=D(A^{\infty})$ and
$x\geq0$. Then by the above formulas%
\[
\lbrack\mu,\sigma]=-\beta\rho\lambda.
\]
The singular set $ \Sigma $ is characterized by%
\[
\frac{d}{dx}h+\frac{\rho^{2}}{\beta}\exp(-\beta x)(1-\exp(-\beta
x))=c\exp(-\beta x)
\]
for some real $c$. Therefore $a_{2}$ is some fixed value, namely $a_{2}%
=\frac{\rho^{2}}{2}$ and $a_{1},a_{3}$ are arbitrary. Consequently the singular $ \Sigma $ set is
an affine subspace for the fixed values $\rho,\beta$:%
\[
h=a_{1}-\frac{\rho^{2}}{2}\Lambda^{2}+a_{3}\Lambda.
\]
Going back to traditional notations for the Vasicek model we write%
\begin{align*}
\Lambda(x)  & =\frac{1}{\beta}(1-\exp(-\beta x))\\ B_{V}(x)  & =\Lambda^{\prime}(x)=e^{-\beta x}\\ A_{V}(x)  &
=b\Lambda(x)-\frac{\rho^{2}}{2}\Lambda(x)^{2},
\end{align*}
then $h$ lies in the singular set $\Sigma$ if and only if%
\[
h\in A_{V}+\left\langle B_{V}\right\rangle
\]
for some value $b$ (which becomes an additional parameter in the short
rate
equation). The solution for $r^{\ast}$ in the singular set reads as follows%
\begin{align*}
r_{t}  & =A_{V}+B_{V}R_{t}\\ dR_{t}  & =(b-\beta R_{t})\,dt+\rho \,dW_{t},
\end{align*}
where $R_{t}=ev_{0}(r_{t})$ denotes the short rate, which is the Vasicek short rate model.

Outside the singular set $ \Sigma $ we have a $2$-dimensional realization. First we calculate the
deterministic part of the dynamics%
\begin{align*}
Fl_{u_{0}}^{\nu}(r^{\ast})(x)  & =S_{u_{0}}r^{\ast}(x)+\int_{0}^{u_{0}}S_{u_{0}%
-s}(\frac{\rho^{2}}{\beta}\exp(-\beta x)(1-\exp(-\beta x))\,ds\\
& =S_{u_{0}}r^{\ast}(x)+\frac{\rho^{2}}{2}\int_{0}^{u_{0}}\frac{d}{dx}%
(\Lambda)^2(x+u_{0}-s)\,ds\\
& =r^{\ast}(x+u_{0})+\frac{\rho^{2}}{2}\Lambda(x+u_{0})^{2}-\frac{\rho^{2}}%
{2}\Lambda(x)^{2}.
\end{align*}
If we identify $u_{0}$ with the time variable $t$, which is possible since the stochastics only occurs in direction
of $B_{HWV}$ (see Remark~\ref{rem1}), we obtain by direct calculations for \eqref{atsm}
\begin{align*}
r_{t}(x)  & =r^{\ast}(x+t)+\frac{\rho^{2}}{2}\Lambda(x+t)^{2}-\frac{\rho^{2}}%
{2}\Lambda(x)^{2}+\Lambda^{\prime}(x)Z_{t}\\ dZ_{t}  & =-\beta
Z_{t}\,dt+\rho \,dW_{t}.
\end{align*}
A parameter transformation yields the customary form, namely $$ R_t = e^{-\beta t}r^{\ast}(0) +
\int_{0}^{t}e^{-\beta(t-s)}b(s)\,ds + Z_t. $$ This yields the following expressions:
\begin{align*}
A_{HWV}(t,x)  & = r^{\ast}(x+t)+\frac{\rho^{2}}{2}\Lambda(x+t)^{2}-\frac{\rho^{2}}%
{2}\Lambda(x)^{2}- \\ &\quad - ({\Lambda}^{\prime}(x))^2 r^{\ast}(0)  - {\Lambda}^{\prime}(x)\int_0^t e^{-\beta(t-s)} b(s)\,ds \\
B_{HWV}(x) & = B_V (x) =\Lambda^{\prime}(x)\\ dR_{t} & =(b(t)-\beta R_{t})\,dt+\rho \,dW_{t}\\ r_{t} &
=A_{HWV}(t)+B_{HWV}R_{t}\\
b(t)  & =\frac{d}{dt}r^{\ast}(t)+\beta r^{\ast}(t)+\frac{\rho^{2}}{2\beta}%
(1-\exp(-2\beta t)).
\end{align*}
The functions $ A_{HWV} $ and $ B_{HWV} $ are solutions of time-dependent Riccati equations constructed by
geometric methods. The equation for $ b $ follows from the fact that $ A_{HWV} (t,0) = 0 $.

\subsection{The Hull-White extension of the CIR model}\label{subsecHWCIR}

We proceed in the same spirit: $\sigma(r):=\rho\sqrt{ev_{0}(r)}\lambda$ for $\rho>0$. The volatility structure is
defined on the convex open set $\Ucal=\{ev_{0}(r)>\epsilon\}$ for some $\epsilon>0$. The function
$\Lambda:=\int\lambda$ satisfies (in certain normalization) a
Riccati equation, namely%
\[
\frac{d}{dx}\Lambda+\frac{\rho^{2}}{2}\Lambda^{2}+\beta\Lambda=1,\quad \Lambda(0)=0.
\]
We obtain the solution (see e.g.~\cite[Section~7.4.1]{fil:lnm01})%
\[
\Lambda(x)=\frac{2\exp(x\sqrt{\beta^{2}+2\rho^{2}})-1}{(\sqrt{\beta^{2}%
+2\rho^{2}}-\beta)(\exp(x\sqrt{\beta^{2}+2\rho^{2}})-1)+2\sqrt{\beta^{2}%
+2\rho^{2}}}.
\]
Under this assumption we can proceed as above: The singular set $ \Sigma $ is determined
by the equation%
\[ \nu(h)=\frac{d}{dx}h+\rho^{2}ev_{0}(h)\Lambda\Lambda^{\prime}\in\left\langle \lambda\right\rangle\]
hence
\[
h=a_{1}+\frac{\rho^{2}}{2}a_{1}\Lambda^{2}+a_{3}\Lambda.\]
Again $a_{1}$ and $a_{3}$ can be chosen freely, which completely determines $\Sigma$. Traditionally one writes the singular set in the following form:%
\begin{align*}
A_{CIR}  & =b\Lambda\\
B_{CIR}  & =1-\beta\Lambda-\frac{\rho^{2}}{2}\Lambda^{2}=\Lambda^{\prime}%
\end{align*}
with some additional parameter $b$ and we obtain equally that $h$ lies
in $\Sigma$ if and only if%
\[
h\in A_{CIR}+\left\langle B_{CIR}\right\rangle .
\]
The short rate dynamics follows the known pattern:%
\begin{align*}
r_{t}  & =A_{CIR}+B_{CIR}R_{t}\\
dR_{t}  & =(b-\beta R_{t})\,dt+\rho\sqrt{R_{t}}\,dW_{t}%
\end{align*}
for $r^{\ast}\in\Sigma$. Outside the singular set we have a
$2$-dimensional realization. First we calculate the deterministic part, by the variation of constants formula,
\[
Fl_{u_{0}}^{\nu}(r^{\ast})(x)=S_{u_{0}}r^{\ast}(x)+\rho^{2}\int_{0}^{u_{0}}%
Fl_{s}^{\nu}(r^{\ast})(0)(S_{u_{0}-s}(\Lambda^{\prime}\Lambda))(x)\,ds.
\]
Identifying $u_{0}$ with the time parameter yields the following formula $2$-dimensional realization, which is derived by direct calculations,
\begin{align*}
r_{t}  & =Fl_{t}^{\nu}(r^{\ast})+\Lambda^{\prime}Z_{t}\\ dZ_{t}  & =-\beta Z_{t}\,dt+\rho\sqrt{c(t)+Z_{t}}\,dW_{t},
\end{align*}
where $c(t)=Fl_{t}^{\nu}(r^{\ast})(0)$. The short rate is given through
$R_{t}=c(t)+Z_{t}$ and%
\begin{align*}
dR_{t}  & =(\beta c^{\prime}(t)-\beta Z_{t})\,dt+\rho\sqrt{R_{t}}\,dW_{t}\\ & =(b(t)-\beta
R_{t})\,dt+\rho\sqrt{R_{t}}\,dW_{t}.
\end{align*}
Notice that $ \lambda(0) = {\Lambda}^{\prime}(0)=1 $ by the Riccati equation and $ b(t) = c'(t) + \beta c(t) $.

This formula closes the circle with the classical Hull-White extension of the
CIR-model:%
\begin{align*}
A_{HWCIR}(t,x) & = Fl_{t}^{\nu}(r^{\ast})(x) - c(t)\Lambda^{\prime}(x)   \\ B_{HWCIR}  & =B_{CIR} = \Lambda^{\prime}\\
dR_t & =  (b(t)-\beta R_{t})\,dt+\rho\sqrt{R_{t}}\,dW_{t} \\ r_{t}  & =A_{HWCIR}(t)+B_{HWCIR}R_{t}\\ b(t) & =\beta c(t)
+ \frac{d}{dt}c(t) \\ c(t) & = r^{\ast}(t) + {\rho}^2 \int_0^t c(s) (\Lambda {\Lambda}^{\prime})(t-s)\,ds.
\end{align*}
Again this is a geometrical construction of solutions of time-dependent Riccati equations.

\subsection{Fitting procedures as leaves of foliations}

A popular forward curve-fitting method is the
Svensson~\cite{sve:94} family
\[
G_{S}(x,z)=z_{1}+z_{2}e^{-z_{5}x}+z_{3}xe^{-z_{5}x}+z_{4}xe^{-z_{6}x}.
\]
It is shown in~\cite{fil:98} that the only non-trivial interest
rate model that is consistent with the Svensson family is of the
form
\begin{equation}
r_{t}=Z_{t}^{1}g_{1}+\dots+Z_{t}^{4}g_{4}, \label{sven}%
\end{equation}
where
\[
g_{1}(x)\equiv1,\quad g_{2}(x)=e^{-\alpha x},\quad
g_{3}(x)=xe^{-\alpha x},\quad g_{4}(x)=xe^{-2\alpha x},
\]
for some fixed $\alpha>0$. Moreover,
\[
Z_{t}^{1}\equiv Z_{0}^{1},\quad Z_{t}^{3}=Z_{0}^{3}e^{-\alpha
t},\quad Z_{t}^{4}=Z_{0}^{4}e^{-2\alpha t}\quad(Z_{0}^{4}\geq0)
\]
and $Z^{2}$ satisfies
\begin{equation}
dZ_{t}^{2}=\left(  Z_{t}^{3}+Z_{t}^{4}-\alpha Z_{t}^{2}\right)
\,dt+\sqrt
{\alpha Z_{t}^{4}}\,dW_{t}. \label{sven2}%
\end{equation}
Here $W$ is a real-valued Brownian motion.

We now shall find a 2-dimensional local HJM model that is
of the form (\ref{sven}) whenever $r_{0}=\sum_{j=1}^{4} z_{j}
g_{j}$ with $z_{4}\ge0$. In view of (\ref{sven2}), a candidate for
$\sigma$ is given on the half space $\Ucal:=\{\ell>0\}$ by
\[
\sigma(h)=\sqrt{\alpha\ell(h)} g_{2},
\]
where $\ell$ is some continuous linear functional on $C(\mathbb{R}_{\ge0},\mathbb{R})$ with
$\ell(g_{1})=\ell(g_{2})=\ell (g_{3})=0$ and $\ell(g_{4})=1$.
Straightforward calculations show, for $h\in \Ucal\cap D(A^\infty)$,
\begin{align*}
\mu(h)  &  =Ah+\ell(h)g_{2}-\ell(h)g_{2}^{2}\\ [\mu,\sigma](h)  &
=-\alpha\sqrt{\alpha\ell(h)}g_{2}-\frac{\ell(\mu
(h))}{2\sqrt{\alpha\ell(h)}}g_{2}.
\end{align*}
(the clue is that $\ell\circ\sigma\equiv0$). Hence indeed the Lie algebra generated by $\sigma$ and $\mu$ has dimension 2 on $\Ucal\cap D(A^\infty)\setminus\Sigma$.

%


\providecommand{\bysame}{\leavevmode\hbox to3em{\hrulefill}\thinspace}
\providecommand{\MR}{\relax\ifhmode\unskip\space\fi MR }
\providecommand{\MRhref}[2]{%
  \href{http://www.ams.org/mathscinet-getitem?mr=#1}{#2}
}
\providecommand{\href}[2]{#2}

\end{document}